\documentclass{amsart}
\usepackage{amssymb}
\usepackage{amscd}
\usepackage{amsmath}
\usepackage{epsfig}
\usepackage{graphicx}

\usepackage[all]{xy}

\usepackage{pgf,tikz}
\usepackage{caption}
\usepackage{shuffle}
 
\usetikzlibrary{arrows}
\usetikzlibrary{decorations.pathreplacing}
\title[Quantum dilogarithm and associators]{Quantum dilogarithm and associators}

\author{Hidekazu Furusho}
\address{Graduate School of Mathematics, Nagoya University, Chikusa-ku, Furo-cho, Nagoya, 464-8602,  Japan}
\email{furusho@math.nagoya-u.ac.jp}


\date{November 30, 2026}

\newcommand{\C}{\mathbb C}
\newcommand{\ii}{\mathrm i}
\newcommand{\ee}{\mathrm e}
\newcommand{\Li}{\operatorname{Li}}
\newcommand{\Assoc}{\operatorname{Assoc}}
\newcommand{\KZ}{\mathrm{KZ}}

\newcommand{\tlie}{\mathfrak t}

\theoremstyle{plain}
\newtheorem{theorem}{Theorem}%[section]
\newtheorem{proposition}[theorem]{Proposition}
\newtheorem{lemma}[theorem]{Lemma}

\theoremstyle{definition}
\newtheorem{remark}[theorem]{Remark}
\numberwithin{equation}{section}
            
\begin{document}
\bibliographystyle{amsalpha+}

%%%%%%%%%%%%%%%%%%%%%%%%%%%%%%%%%%%%%%%%%%%%%%%%%%%%%%%%%%%%%%%%%%%%%%
\begin{abstract}
We prove an explicit specialization formula expressing the quantum dilogarithm in terms of any associator. 
One generator acts by a multiple of the Euler derivation,  and the other acts by its negative plus a logarithmic multiplier.
The resulting operator is independent of the choice of associator and equals multiplication by the quantum dilogarithm times an exponential of the classical dilogarithm.
The proof is based on a comparison of the two 
hexagon equations.
\end{abstract}
\maketitle

%\tableofcontents

%\setcounter{section}{-1}
\section{Introduction}\label{sec:intro}

Pentagon identities arise in three settings relevant to this
paper: the classical five-term identity in Rogers form of the classical dilogarithm, the
noncommutative pentagon identity for the quantum dilogarithm
\cite{FK,KN}, and Drinfeld's pentagon equation for associators
\cite{Dr}. The classical five-term identity arises as a
semiclassical limit of the quantum pentagon \cite{KN}.
These parallels motivate the
question we address here: can the quantum dilogarithm be
recovered by specializing an associator?

We evaluate an associator at two first-order differential
operators whose sum acts by multiplication by a scalar
multiple of \(\Li_1(-X)=-\log(1+X)\). For a fixed nonzero
parameter \(\mu\), the resulting operator is independent of
the choice of associator. It acts by multiplication by the
quantum dilogarithm times the exponential of a scalar
multiple of \(\Li_2(-X)\). We prove the formula by comparing
the hexagon relation at \(\mu\) and at \(-\mu\) for the
same associator.

Let \(\tau\) be a formal variable, set \(K=\C((\tau))\) and
\(R=K[[X]]\), and put \(q_\mu=\ee^{\mu\tau/2}\) for
\(\mu\in\C^\times\). We use the following normalization of
the compact quantum dilogarithm:
\begin{equation}\label{eq:series}
 \Psi_q(X)=\sum_{n\geq0}
       \frac{(-q)^n X^n}{\prod_{r=1}^{n}(1-q^{2r})},
\end{equation}
where the empty product is \(1\). For \(q=q_\mu\), this
defines a series in \(R\). For complex \(q\) with \(|q|<1\),
the series is the Taylor expansion at \(X=0\) of
\begin{equation}\label{eq:product}
 \prod_{k\geq0}(1+q^{2k+1}X)^{-1}.
\end{equation}

We write \(\Assoc_\mu(\C)\) for the set of associators over
\(\C\) with parameter \(\mu\); see
Section~\ref{sec:associator} for the defining relations.
Set
\begin{equation}\label{eq:introdefs}
 D=X\partial_X,\qquad
 \Li_s(-X)=\sum_{m\geq1}\frac{(-1)^mX^m}{m^s}
 \quad (s=1,2).
\end{equation}
We denote multiplication by \(g\) on \(R\) by \(M_g\),
also written \(M[g]\).

\begin{theorem}\label{thm:main}
For every \(\mu\in\C^\times\) and every
\(\Phi\in\Assoc_\mu(\C)\), we have
\begin{equation}\label{eq:generalmu}
 \Phi\!\left(\tau D,-\tau D+\frac{1}{\mu}M_{\Li_1(-X)}\right)
 =M\!\left[
     \exp\!\left(\frac{\Li_2(-X)}{\mu\tau}\right)
     \Psi_{q_\mu}(X)
   \right].
\end{equation}
The two factors on the right are defined in \(K[[X]]\),
and their product lies in \(\C[[\tau]][[X]]\).
\end{theorem}

For the KZ associator, we use the convention
\[
 G'=\left(\frac{\mathbf e_0}{z}
          +\frac{\mathbf e_1}{z-1}\right)G,\qquad
 \Phi_{\KZ}=G_1^{-1}G_0,
\]
where \(G_0\sim z^{\mathbf e_0}\) at \(0\) and
\(G_1\sim(1-z)^{\mathbf e_1}\) at \(1\).
With this convention,
\(\Phi_{\KZ}\in\Assoc_{2\pi\ii}(\C)\) \cite{Dr}.
Thus, setting \(q=\ee^{\pi\ii\tau}\), we obtain
\begin{equation}\label{eq:main}
 \Phi_{\KZ}\!\left(
     \tau D,-\tau D+\frac{1}{2\pi\ii}M_{\Li_1(-X)}
   \right)
 =M\!\left[
     \exp\!\left(\frac{\Li_2(-X)}{2\pi\ii\tau}\right)
     \Psi_q(X)
   \right].
\end{equation}

%%%%%%%%%%%%%%%%%%%%%%%%%
\section{Associators and the hexagon equations}\label{sec:associator}
Let \(\C\langle\!\langle\mathbf e_0,\mathbf e_1
\rangle\!\rangle\) be the degree-completed free associative algebra.
Let \(\tlie_4\) be the Lie algebra generated by
\(t_{ij}=t_{ji}\), subject to
\[
 [t_{ij},t_{ik}+t_{jk}]=0,\qquad [t_{ij},t_{kl}]=0
\]
for pairwise distinct indices in each relation.

Following \cite[\S5]{Dr}, for \(\mu\in\C^\times\), let
\(\Assoc_\mu(\C)\) denote the set of group-like series
\(\Phi=\Phi(\mathbf e_0,\mathbf e_1)\) with constant term \(1\)
satisfying the pentagon equation
\begin{equation}\label{eq:drpent}
 \begin{aligned}
 &\Phi(t_{12},t_{23}+t_{24})\Phi(t_{13}+t_{23},t_{34})
 %\\
 %&\qquad
 =\Phi(t_{23},t_{34})
 \Phi(t_{12}+t_{13},t_{24}+t_{34})\Phi(t_{12},t_{23})
 \end{aligned}
\end{equation}
in \(\widehat U(\tlie_4)\), the duality relation
\begin{equation}\label{eq:duality}
 \Phi(\mathbf e_0,\mathbf e_1)
 \Phi(\mathbf e_1,\mathbf e_0)=1,
\end{equation}
and the hexagon equation
\begin{equation}\label{eq:threecycle}
 \Phi(\mathbf e_0,\mathbf e_1)\ee^{\mu\mathbf e_0/2}
 \Phi(\mathbf e_\infty,\mathbf e_0)\ee^{\mu\mathbf e_\infty/2}
 \Phi(\mathbf e_1,\mathbf e_\infty)\ee^{\mu\mathbf e_1/2}=1,
\end{equation}
where \(\mathbf e_\infty=-\mathbf e_0-\mathbf e_1\).

In \cite{Fu}, it is shown that 
the pentagon equation  \eqref{eq:drpent} for a group-like series  \(\Phi\)
forces the linear
terms of \(\Phi\) to vanish and implies  \eqref{eq:duality}
 and \eqref{eq:threecycle}
 for a parameter \(\mu\) satisfying
\begin{equation}\label{eq:mu2}
 \mu^2=24\langle\Phi\mid\mathbf e_0\mathbf e_1\rangle.
\end{equation}
Here \(\langle\Phi\mid w\rangle\) denotes the coefficient
of the word \(w\) in \(\Phi\).

Interchanging \(\mathbf e_0\) and \(\mathbf e_1\) in
\eqref{eq:threecycle}, inverting the product, and using
\eqref{eq:duality}, followed by a cyclic rearrangement,
gives the same relation with \(\mu\) replaced by \(-\mu\).
Thus, for the same \(\Phi\), we have the signed hexagon equations
\begin{equation}\label{eq:hexagon}
 \begin{gathered}
 \Phi(\mathbf e_0,\mathbf e_1)\ee^{\varepsilon\mu\mathbf e_0/2}
 \Phi(\mathbf e_\infty,\mathbf e_0)
 \ee^{\varepsilon\mu\mathbf e_\infty/2}
 \Phi(\mathbf e_1,\mathbf e_\infty)
 \ee^{\varepsilon\mu\mathbf e_1/2}=1,
 \qquad \varepsilon\in\{1,-1\}.
 \end{gathered}
\end{equation}

%%%%%%%%%%%%%%%%%%%%%%%%%%%%%%%
\section{Proof of the specialization theorem}\label{sec:hex}
We first treat an arbitrary multiplier \(\ell\in X\C[[X]]\).
Fix \(\mu\in\C^\times\) and \(\Phi\in\Assoc_\mu(\C)\),
and write \(q=q_\mu=\ee^{\mu\tau/2}\).
%It is useful to leave the multiplier unspecified at first. Fix \(\mu\in\C^\times\), let \(q=\ee^{\mu\tau/2}\), and take any \(\ell(X)\in X\C[[X]]\). 
Set
\begin{equation}\label{eq:rep}
 A=\tau D,\qquad B=-\tau D-M_\ell,\qquad C=M_\ell .
\end{equation}
We first regard \(A,B,C\) as \(\C[[\tau]]\)-linear
endomorphisms of \(R_0=\C[[\tau,X]]\). Since they raise the
\((\tau,X)\)-adic order, we can evaluate \(\Phi\) at pairs
of these operators on \(R_0\). The resulting endomorphisms
preserve the \(X\)-adic filtration and therefore extend to
\(R=K[[X]]\) after inverting \(\tau\) and completing
\(X\)-adically. All subsequent operator identities are
taken on \(R\).

\begin{lemma}\label{lem:mult}
Each of \(\Phi(A,B),\Phi(C,A),\Phi(B,C)\) acts  by multiplication by an element of \(1+XR_0\).
\end{lemma}
\begin{proof}
We have \([A,M_g]=M_{\tau Dg}\) and \([M_g,M_h]=0\). Thus every Lie bracket of length at least two in \eqref{eq:rep} is multiplication by an element of \(XR_0\). 
Since \(\Phi\) is group-like and has no linear term,
\(\log\Phi\) is a Lie series of degree at least two.
This applies to the image of \(\log\Phi\) for each of the three pairs. Exponentiating proves the assertion.
\end{proof}

Write \(\Phi(A,B)=M_{p_\ell}\) with  $p_\ell\in 1+XR_0$. We seek a multiplication
operator \(J=M_j\), with \(j\in 1+XR\), such that
\begin{equation}\label{eq:conjugacy}
 B=J(-A)J^{-1}.
\end{equation}
By the Leibniz rule, for \(h\in R\),
\[
 J(-A)J^{-1}h
 =-\tau jD(j^{-1}h)
 =-\tau Dh+\tau(D\log j)h.
\]
Since \(B=-A-M_\ell\), equation \eqref{eq:conjugacy}
is equivalent to \(\tau D\log j=-\ell\).
On \(XR\), the operator \(D\) is invertible, with
\(D^{-1}X^m=X^m/m\) for \(m\geq 1\). Thus the unique
solution with \(j(0)=1\) is \(j_\ell\) as defined below.
With \(J=M_{j_\ell}\), we set
\begin{equation}\label{eq:gauge}
 j_\ell:=\exp(-\tau^{-1}D^{-1}\ell),
 \qquad
 F_\ell:=j_\ell^{-1}p_\ell.
\end{equation}
Then \(F_\ell\in 1+XR\) and
\(M_{F_\ell}=J^{-1}\Phi(A,B)\).

\begin{proposition}\label{prop:eliminate}
Let $\mu\in\C^\times$ and $q=q_\mu$.
For every \(\Phi\in\Assoc_\mu(\C)\), the multiplier \(F_\ell\) is the unique series in \(1+XK[[X]]\) satisfying
\begin{equation}\label{eq:centered}
 \frac{F_\ell(qX)}{F_\ell(q^{-1}X)}=\ee^{\mu\ell(X)} .
\end{equation}
In particular, if \(\ell=\sum_{m\ge1}\ell_mX^m\), then
\begin{equation}\label{eq:kernel}
 \log p_\ell(X)=\sum_{m\ge1}\ell_m\,H_\mu(m\tau)X^m,
 \qquad H_\mu(s)=\frac{\mu}{2\sinh(\mu s/2)}-\frac1s .
\end{equation}
\end{proposition}

\begin{proof}
Put \(T_\varepsilon=\ee^{\varepsilon\mu A/2}\), so that \(T_\varepsilon h(X)=h(q^\varepsilon X)\). Conjugation gives
\begin{equation}\label{eq: conjugation formulas}
 T_\varepsilon M_gT_\varepsilon^{-1}=M_{g(q^\varepsilon X)},\qquad
 \ee^{\varepsilon\mu B/2}=JT_\varepsilon^{-1}J^{-1}.
\end{equation}
The first identity follows by applying both sides to \(h\in R\),
since \(T_\varepsilon h(X)=h(q^\varepsilon X)\).
For the second, work modulo \(X^{N+1}\). On this quotient,
the exponential series for \(-\varepsilon\mu A/2\) converges
\(\tau\)-adically, and conjugation by \(J\) is continuous in
the space of matrices over \(K\). Thus \(B=J(-A)J^{-1}\)
implies
\(\ee^{\varepsilon\mu B/2}
=J\ee^{-\varepsilon\mu A/2}J^{-1}
=JT_\varepsilon^{-1}J^{-1}\)
on each quotient, and hence on \(R\).

By Lemma~\ref{lem:mult}, we may define \(k\in1+XK[[X]]\) by
\[
 M_k=\Phi(C,A)\Phi(B,C)J.
\]
Substituting the conjugation formulas \eqref{eq: conjugation formulas} into
\eqref{eq:hexagon}, conjugating by  \(J\), and using
\(J^{-1}\Phi(A,B)=M_{F_\ell}\), we obtain
\[
 M_{F_\ell}T_\varepsilon
 \bigl(\Phi(C,A)e^{\varepsilon\mu C/2}\Phi(B,C)J\bigr)
 T_\varepsilon^{-1}=1.
\]
The factors inside parentheses are multiplication operators
and hence commute. By the definition of \(k\) and the equality
\(C=M_\ell\), their product is
\[
 \Phi(C,A)\Phi(B,C)J\,e^{\varepsilon\mu C/2}
 =M_kM_{\exp(\varepsilon\mu\ell/2)}
 =M_{k\exp(\varepsilon\mu\ell/2)}.
\]
Thus
\[
 M_{F_\ell}\,T_\varepsilon\,
 M_{k\exp(\varepsilon\mu\ell/2)}
 \,T_\varepsilon^{-1}=1.
\]
By the first identity in \eqref{eq: conjugation formulas},
this operator identity is equivalent to
\[
 F_\ell(X)\,k(q^\varepsilon X)\,
 \exp\!\left(\frac{\varepsilon\mu}{2}\ell(q^\varepsilon X)\right)=1.
\]
In the positive identity replace \(X\) by \(q^{-1}X\), and in the negative identity replace \(X\) by \(qX\). The two equations now read
\begin{align*}
 F_\ell(q^{-1}X)\,k(X)\,\ee^{\mu\ell(X)/2}&=1,\\
 F_\ell(qX)\,k(X)\,\ee^{-\mu\ell(X)/2}&=1 .
\end{align*}
The unknown \(k\) occurs at the same argument. Taking the quotient gives \eqref{eq:centered}.

Taking formal logarithms in \eqref{eq:centered} yields
\begin{equation}\label{eq: coefficient formula}
 (q^m-q^{-m})\langle \log F_\ell \mid  X^m\rangle=\mu\ell_m \qquad (m\ge1).
\end{equation}
Here \(\langle f\mid X^m\rangle\) denotes the coefficient of \(X^m\) in \(f\).
The factor on the left is nonzero in \(K\), which proves uniqueness. Adding 
\(\langle \log j_\ell  \mid X^m\rangle =-\ell_m/(m\tau)\) gives \eqref{eq:kernel}.
\end{proof}

\begin{remark}
The specialization \eqref{eq:rep} factors through the
metabelian quotient, since its derived algebra consists of
commuting multiplication operators. Write the class of
\(\log\Phi\) as
\(f_\Phi(\operatorname{ad}\mathbf e_0,\operatorname{ad}\mathbf e_1)
[\mathbf e_0,\mathbf e_1]\).
For \(\ell=X\), we have \([A,B]=-\tau M_X\), and
\(\operatorname{ad}A\) and \(\operatorname{ad}B\) act on \(M_X\)
by \(\tau\) and \(-\tau\), respectively. 
%Thus \eqref{eq:kernel} gives \(-\tau f_\Phi(\tau,-\tau)=H_\mu(\tau)\).
For \(\ell=X\), only the \(m=1\) term in \eqref{eq:kernel}
remains, giving \(\log\Phi(A,B)=H_\mu(\tau)M_X\).
Comparing the two expressions for \(\log\Phi(A,B)\) gives
\(-\tau f_\Phi(\tau,-\tau)=H_\mu(\tau)\).
The two hexagons therefore determine \(f_\Phi(s,-s)\) for fixed
\(\mu\) without using the pentagon.
Since \(H_\mu(s)\) is odd, \(f_\Phi(s,-s)=\sum_{j\ge0}c_j s^{2j}\).
Thus
\[
\log\Phi(A,B)
=\sum_{j\ge0}c_j(\operatorname{ad}A)^{2j}[A,B]
=\sum_{j\ge0}c_j(\operatorname{ad}A)^{2j+1}B.
\]
For the KZ associator \(\Phi_{\mathrm{KZ}}\),
Drinfeld's formula
\cite[Eq.~(2.15)]{Dr} gives
\[
 1-s^2 f_{\Phi_{\mathrm{KZ}}}(s,-s)
 =\Gamma(1-s)\Gamma(1+s)
 =\frac{\pi s}{\sin(\pi s)}
 =1+sH_{2\pi\ii}(s).
\]
\end{remark}

\begin{proof}[Proof of Theorem~\ref{thm:main}]
%We first allow arbitrary \(\mu\ne0\) and 
Take 
\[\ell(X)=\mu^{-1}\log(1+X)=-\mu^{-1}\Li_1(-X).\] 
Since \(D\Li_2(-X)=\Li_1(-X)\), \eqref{eq:gauge} gives 
\[j_\ell=\exp(\Li_2(-X)/(\mu\tau)).\] 
%Write \(F=F_\ell\). 
Replacing \(X\) by \(qX\) in \eqref{eq:centered} gives
\begin{equation}\label{eq:difference}
 F_\ell(q^2X)=(1+qX)F_\ell(X),\qquad F_\ell(0)=1 .
\end{equation}
If \(F_\ell=\sum_{n\ge0}b_n X^n\), then \((q^{2n}-1)b_n=qb_{n-1}\) for \(n\ge1\). 
Thus \(b_n=(-q)^n/\prod_{r=1}^n(1-q^{2r})\) and
 \[F_\ell=\Psi_q.\] 
 Since \(p_\ell=j_\ell F_\ell\), this proves \eqref{eq:generalmu}.
By \eqref{eq:kernel} and \(H_\mu(s)=-\mu^2s/24+O(s^3)\),
we have \(\log p_\ell\in\tau X\C[[\tau]][[X]]\). In particular,
\[
 \log\Psi_q(X)
 =\log p_\ell(X)-\frac{\Li_2(-X)}{\mu\tau}
 =-\frac{\Li_2(-X)}{\mu\tau}+O(\tau)
 \qquad\text{coefficientwise in }X.
\]
Thus the exponential factor in \eqref{eq:generalmu} cancels
the pole in the logarithm, and the product belongs to
\(\C[[\tau]][[X]]\). Taking \(\mu=2\pi\ii\) gives
\eqref{eq:main}.
\end{proof}

\begin{remark}
The proof uses the two cases of \eqref{eq:hexagon},
together with the fact that \(\log\Phi\) is a Lie series
with no linear term. The pentagon equation
\eqref{eq:drpent} is not needed directly: for a group-like
solution of the pentagon equation, the result in 
\cite{Fu} supplies the relations used here, with the
parameter determined up to sign by \eqref{eq:mu2}.

Together with Sch\"utzenberger's addition formula, the difference
equation \eqref{eq:difference} yields the standard pentagon
identity for the quantum dilogarithm; see \cite{Vo,KN}.
This gives an indirect route from Drinfeld's pentagon to the
quantum dilogarithm pentagon, via the hexagon relations and
the \(q\)-difference equation \eqref{eq:difference}.
\end{remark}

\medskip
{\bf Declaration of Generative AI Use.}
This work grew out of the author's exploratory discussions with ChatGPT-6.2-Astra (OpenAI); 
in particular, the specialization \eqref{eq:rep} was suggested by the model. 
It also assisted with the exposition  of the manuscript.
All mathematical content was independently checked by the author, who bears sole responsibility for this article.in particular, the specialization \eqref{eq:rep} was suggested by the model. 

\bigskip
{\bf Acknowledgements.}
 This work is supported by the following grants: JSPS KAKENHI JP24K00520, and  JP24K21510. 

%%%%%%%%%%%%%%%%%%%%%%%%%%%%%

\end{document}